\documentclass[11pt]{article}
\usepackage[hmargin=1.04in,top=0.92in,bottom=0.92in]{geometry}
\usepackage{amsmath,amssymb,amsthm}
\usepackage[hidelinks]{hyperref}
\hypersetup{pdftitle={The perfect set property for ordinal-definable sets in single-real extensions}}
\numberwithin{equation}{section}
\newtheorem{theorem}{Theorem}[section]
\newtheorem{proposition}[theorem]{Proposition}
\newtheorem{lemma}[theorem]{Lemma}
\newtheorem{corollary}[theorem]{Corollary}
\theoremstyle{remark}
\newtheorem*{remark}{Remark}
\newcommand{\PP}{\mathbb P}
\newcommand{\BB}{\mathbb B}
\newcommand{\QQ}{\mathbb Q}
\newcommand{\WT}{\mathbb W}
\DeclareMathOperator{\Coll}{Coll}
\newcommand{\RO}{\operatorname{RO}}
\newcommand{\stem}{\operatorname{stem}}
\newcommand{\restr}{\mathbin{\upharpoonright}}
\newcommand{\forces}{\Vdash}
\newcommand{\PK}{\mathrm{PK}(\mathrm{OD})}
\title{The perfect set property for ordinal-definable sets\\in single-real extensions}
\author{Jason Zesheng Chen\thanks{I thank Toby Meadows for extensive discussions of this problem.}}
\date{}
\begin{document}
\maketitle
\begin{abstract}
Kanovei and Lyubetsky asked whether the perfect set property for ordinal-definable sets of
reals is consistent with all reals belonging to a single $L[a]$, and whether it holds after
the single L\'evy collapse of $\omega_1^L$ to $\omega$. We give a positive
consistency result for the former using Prikry's wide-tree collapse and a negative answer for the
latter. The methods used for the latter also provide a negative answer to a recent question posed by Schilhan.
\end{abstract}

\section{Introduction}

Write $\PK$ for the assertion that every uncountable ordinal-definable subset of
$\omega^\omega$ contains a nonempty perfect subset. Equivalently, it must contain a copy of
$2^\omega$. Kanovei and Lyubetsky ask whether $\PK$ can hold when
$\omega^\omega\subseteq L[a]$ for some real $a$ and, in particular, whether it holds in
the extension of $L$ by the L\'evy collapse of $\omega_1^L$ to $\omega$
\cite[\S7B and Remark~7.14]{KL03}.
They proved that $\PK$ holds in the extension by a product of $(\omega_2)^L$ collapse
forcings \cite[Theorem~7.9]{KL03}.

\begin{theorem}\label{thm:summary}
Let $\kappa=\omega_1^L$. We have the following:
\begin{enumerate}
\item In the extension by the full wide-tree forcing on $\kappa^{<\omega}$, $\PK$ holds.
\item In the extension by the L\'evy collapse $\Coll(\omega,\kappa)$, $\PK$ fails.
\end{enumerate}
\end{theorem}

In both extensions the generic object is a sequence $g:\omega\to\kappa$. For each
$\xi<\kappa$, let $j(\xi)$ be the $L$-least real coding a well-order of type $\xi$ on
a subset of $\omega$. Fix a recursive pairing function $\langle n,m\rangle$ on $\omega$.
For every $h\in\kappa^\omega$, define the real $J(h)\in\omega^\omega$ by
\[
 J(h)(\langle n,m\rangle)=j(h(n))(m) \qquad (n,m<\omega).
\]
With $\kappa$ discrete, $J$ is continuous for the product topology on $\kappa^\omega$.
The $n$th row of $J(h)$ codes a well-order of type $h(n)$, so $J$ is injective and
$L[J(h)]=L[h]$. Write $N_s$ for the basic cylinder of sequences extending a
finite string $s$. A generic branch is a \emph{full generator} of an extension if adjoining
it to the ground model gives that entire extension.

\section{The wide-tree collapse}\label{sec:wide}

In $L$ let $\WT$ consist of the nonempty subtrees of $\kappa^{<\omega}$ in which every node
extends to a node with $\kappa$ immediate successors. 
A subtree is closed under initial
segments, and a smaller tree is a stronger condition. This is Prikry's collapse by wide trees
\cite{CKM84,Miller09}. For a tree $T$ and $s\in T$, put
\[
 T\restr s=\{t\in T:t\subseteq s\text{ or }s\subseteq t\}.
\]
This is again a condition. Let $[T]$ denote the set of branches of $T$ in the universe under
discussion, and let $\stem(T)$ be the longest node common to all its branches. If $G$ is
generic, the union $g$ of the stems of its conditions is cofinal in $\kappa$. Indeed, below
any condition one can put an arbitrarily large ordinal into its stem by restricting beyond a
node with $\kappa$ successors. The same argument allows arbitrarily long stems. For each
$T\in\WT^L$ it is dense to be below $T$ or to have a stem outside $T$. Consequently
\begin{equation}\label{eq:filter}
 G=\{T\in\WT^L:g\in[T]\},\qquad L[G]=L[g].
\end{equation}
We use the following form of Prikry's theorem, whose statement can be found in
\cite[Theorem~24]{Miller09}. We include the argument for completeness but also because equality of the extensions will
be needed for arbitrary ordinal definitions.

\begin{lemma}[Prikry]\label{lem:prikry}
If $G$ is $\WT$-generic over $L$ and $h\in L[G]\cap\kappa^\omega$ has cofinal range in
$\kappa$, then $h$ is itself a $\WT$-generic branch over $L$ and $L[h]=L[G]$.
\end{lemma}
\begin{proof}
Work in $L$. Fix a name $\dot h$ and a condition $p$ forcing it to be cofinal in $\kappa$. We
construct $q\le p$, a wide tree $S$, and an isomorphism between the forcings below $q$ and
$S$ which sends the input generic branch to $\dot h$. The construction can be carried out
below any strengthening of $p$.

Index the construction by $s\in\kappa^{<\omega}$. At $s$ we have a condition $p_s$, a finite
string $v_s$ with $p_s\forces v_s\subset\dot h$, and a node $u_s\in p_s$ with $\kappa$
immediate successors. Begin with $p_\varnothing=p$ and $v_\varnothing=\varnothing$. For each
successor $u_s{}^\frown\langle\xi\rangle$ strengthen the restriction of $p_s$ there to a
condition deciding an initial segment $w_\xi$ of $\dot h$ properly extending $v_s$ and
containing an entry greater than $\xi$. This is possible by cofinality. There are $\kappa$
distinct $w_\xi$, since the entries occurring in any countable family of finite strings have
supremum below $\kappa$. Thin to $\kappa$ distinct strings of one length. In their
finite-height prefix tree some node $r_s\supseteq v_s$ has $\kappa$ immediate successors.
Otherwise that tree would be countable. Keep one $w_\xi$ through each of those successors,
and index the retained choices by $\eta<\kappa$. Take the corresponding strengthened
conditions and strings as $p_{s{}^\frown\langle\eta\rangle}$ and
$v_{s{}^\frown\langle\eta\rangle}$, and choose a node with $\kappa$ immediate successors
above the designated input successor in each child tree.

Let $q$ and $S$ be the downward closures of the $u_s$ and the $r_s$, respectively. The
children indexed at $s$ leave $u_s$ at distinct immediate successors in $q$ and leave $r_s$
at distinct immediate successors in $S$. Between a designated node and its designated child
there is no further split. These trees are wide, $q\le p$, and following the same sequence of
successor indices gives a homeomorphism $F:[q]\to[S]$ coded in $L$. This correspondence also
sends wide subtrees to wide subtrees in both directions. A wide subtree of $q$ is determined
by the designated nodes $u_s$ that it contains, since it is pruned and every intervening
segment has a unique continuation to the next designated node. Taking the downward closure of
the corresponding $r_s$ preserves the property that every node extends to a node with
$\kappa$ successors. The inverse correspondence is the same construction in the other
direction.

For every $n$, the restrictions of $q$ at nodes $u_s$ with $|s|=n$ are predense below $q$.
Each such restriction is contained in $p_s$ and forces $v_s\subset\dot h$. Along a branch the
$v_s$ increase in length, and each $r_s$ lies between $v_s$ and the next
$v_{s{}^\frown\langle\eta\rangle}$. Thus $q\forces\dot h=F(\dot g)$. The forcing isomorphism
sends a generic filter below $q$ to a generic filter below $S$ with branch $F(g)$, and its
inverse recovers the original filter. Such $q$ are dense below $p$, proving the lemma.
\end{proof}

\begin{proposition}\label{thm:wide}
In $L[G]$, every ordinal-definable $X\subseteq\omega^\omega$ containing a nonconstructible
real contains a copy of $2^\omega$.
\end{proposition}
\begin{proof}
Let $x\in X\setminus L$. Write $X=\{z:\varphi(z,\vec\alpha)\}$ in $L[G]$, and choose a name
$\tau$ for $x$ and $p\in G$ with
\[
 p\forces \tau\in\omega^\omega\setminus L\ \wedge\ \varphi(\tau,\vec\alpha).
\]
Below any $q\le p$ there are two conditions deciding incompatible finite initial segments of
$\tau$. Indeed, if all possible finite decisions below $q$ were compatible, the deciding
conditions for each coordinate would yield a unique possible value. The forcing relation in
$L$ would then define a constructible real forced by $q$ to equal $\tau$.

Work now in $L[G]$, with $g$ as in~\eqref{eq:filter}. Construct ground conditions $p_s$ and
ground strings $w_s\in\omega^{<\omega}$ for $s\in 2^{<\omega}$, starting at
$(p,\varnothing)$. Given $p_s$ at level $n$, strengthen it by extending its stem to include,
at a new position, an ordinal greater than $g(n)$. Below that condition take two refinements
deciding incompatible extensions $w_{s{}^\frown i}$ of $w_s$ in $\tau$, of length at least
$n+1$. Use these refinements as the two child conditions. Each condition belongs to $L$,
while the whole system is constructed in the extension.

For $z\in 2^\omega$ in the extension put
\[
 h_z=\bigcup_n\stem(p_{z\restr n}),\qquad
 y_z=\bigcup_n w_{z\restr n}.
\]
The first sequence is cofinal in $\kappa$ and belongs to the body of every $p_{z\restr n}$.
By Lemma~\ref{lem:prikry} and~\eqref{eq:filter}, its generic filter $G_z$ contains all these
conditions and $L[G_z]=L[G]$. Hence $y_z=\tau[G_z]$ and $L[G]\models\varphi(y_z,\vec\alpha)$.
The incompatible sibling decisions make $z\mapsto y_z$ injective, and the length requirements
make it continuous. Its compact image is a Cantor subset of $X$.
\end{proof}

Every ordinal below $\kappa$ was countable in $L$, so the cofinal sequence $g$ makes $\kappa$
countable. By $\mathsf{CH}$ in $L$, the constructible reals are countable in $L[G]$. Thus
every uncountable OD set of reals satisfies the hypothesis of Proposition~\ref{thm:wide}.
Equation~\eqref{eq:filter} gives $L[G]=L[J(g)]$. This proves the first part of
Theorem~\ref{thm:summary}.

\section{Generators of the L\'evy collapse}\label{sec:levy}

We present $\PP=\Coll(\omega,\kappa)$ as $\kappa^{<\omega}$ ordered by end extension.
We identify its generic filter with the union of its conditions. Every $\PP$-generic
$x$ over $L$ maps $\omega$ onto $\kappa$.

\begin{theorem}\label{thm:single-collapse}
If $g$ is $\PP$-generic over $L$, then in $L[g]$ the set
\[
 \{h\in\kappa^\omega:h\text{ is $\PP$-generic over }L
                         \text{ and }L[h]=L[g]\}
\]
contains no homeomorphic copy of $2^\omega$.
\end{theorem}

We prove Theorem~\ref{thm:single-collapse} in the next two subsections. The argument is motivated by the analogy with cone avoidance in computability theory. Given the collapse generic $x$ and a tree $T_x$ in $L[x]$, we seek a branch $y$ through $T_x$ from which $x$ cannot be recovered, in the sense that $x\notin L[y]$. The key desired property is that overlapping output supports prevent local recovery of the input. This already appeared in a special case in an earlier note on countable OD sets in a Cohen extension \cite{Chen26}. There, two distinct slices of a Boolean product can have overlapping support in the base algebra, in which case the base generic does not determine the slice. Here we arrange such an overlap below every input condition. The main additional step is to obtain the overlap while requiring both local isomorphisms to select branches through $T_x$, and to preserve it as further branch values are decided.

So suppose toward a
contradiction that the full generators contain a copy of $2^\omega$ in $L[g]$.
There are a condition $p_*$ and a name $\dot T$ such that
$p_*$ forces $\dot T\subseteq\kappa^{<\omega}$ to be a pruned finitely branching
tree with nonempty perfect body, all of whose branches are generics generating
the extension. Indeed the tree of initial segments of a copy of Cantor space has
these properties by compactness. Write $x$ for the generic below $p_*$ and
$T_x$ for the interpretation of $\dot T$.

In Subsection~\ref{subsec:local-isomorphisms} we obtain local isomorphisms that select
branches of $T_x$ on disjoint input domains and have a common target cylinder as their range.
In Subsection~\ref{subsec:selecting-a-branch} we use them to construct a single branch
name $f(x)$ such that below every input condition two incompatible refinements have
overlapping target supports. A name recovering $x$ from $f(x)$ on that condition would
force the supports to be disjoint.

\subsection{Local isomorphisms}\label{subsec:local-isomorphisms}

Work in $L$ and put $\BB=\RO(\PP)$. We identify a cylinder with its Boolean
value. Ground open sets in outer models are interpreted from their codes by finite
strings. A name for a $\PP$-generic point in a forcing with Boolean completion
$\mathbb A$ induces a complete homomorphism $i:\BB\to\mathbb A$. Write
\[
 \pi_i(a)=\bigwedge\{b\in\BB:a\le i(b)\}.
\]
Thus $b\wedge\pi_i(a)>0$ precisely when $a\wedge i(b)>0$. The definition also applies when
$i$ has a nonzero kernel.

If $v(x)$ is forced to be generic and a condition forces $x\in L[v(x)]$,
a stronger condition $e$ and a ground target name $H$ satisfy
$e\forces x=H[v(x)]$. The induced homomorphism
$b\mapsto e\wedge i_v(b)$ from $\BB$ to $\BB\mathbin{\upharpoonright}e$ is surjective, since
the inverse name reads every basic input condition from the target. After discarding its
kernel we
obtain an isomorphism
\[
 \BB\mathbin{\upharpoonright}e\ \cong\
 \BB\mathbin{\upharpoonright}\pi_{i_v}(e).
\]
We write $\varphi:E\to D$ for such a local isomorphism, using the direction
from source conditions to target conditions. The same symbol denotes the
induced map on generic points. A total name for a point of $\kappa^\omega$ is continuous on
the dense $G_\delta$ obtained by intersecting the dense open sets deciding its individual
coordinates.

Consider the iteration $\QQ$ which first adds $x$ below $p_*$ and then adds
a Cohen branch $y$ through $T_x$, using its nonempty relatively open subsets.
Basic pair conditions have the form $r=(s,t)$, with $s\le p_*$ and
$s\forces t\in\dot T$. They specify $x\in N_s$ and $y\in N_t$.
The iteration forces $y$ to be $\PP$-generic over $L$. Indeed, for each ground
dense $D\subseteq\PP$, the cylinders from $D$ cover $[T_x]$ in $L[x]$.
A finite subcover still covers all branches of $T_x$ in a further extension, since the
finitely branching tree avoiding that finite subcover has bounded height.
Let $i_y:\BB\to\RO(\QQ)$ be the homomorphism induced by $y$ and abbreviate
$\pi_{i_y}$ to $\pi_y$.

\begin{lemma}\label{lem:levy-pair}
$\QQ$ forces $x\notin L[y]$.
\end{lemma}
\begin{proof}
Otherwise a basic pair condition $(s,t)$ forces $x=F[y]$ for a ground target
name $F$, taken to be total with values in $\kappa^\omega$ by assigning a
default value if necessary. Fix a generic $x\in N_s$. In $L[x]$ we have
$F[z]=x$ for every $z\in[T_x]\cap N_t$. Otherwise, since $z$ is
$\PP$-generic, some finite prefix of $z$ decides a coordinate of $F[z]$
differently from the corresponding coordinate of $x$. Forcing a relative
Cohen branch through this prefix contradicts $(s,t)$.

Fix such a $z$ and choose a ground input name $H$ with $H[x]=z$. In $L[z]$,
$F[z]$ is $\PP$-generic over $L$ and $H[F[z]]=z$. By the truth lemma some
prefix $d$ of $z$ forces both assertions about the canonical target generic.
Thus every generic $z'\in N_d$ with $F[z']=x$ satisfies $z'=H[x]=z$.
This isolates $z$ in $[T_x]\cap N_t$, a contradiction.
\end{proof}

Say a local isomorphism $\varphi:E\to D$ is \emph{through} $r=(s,t)$ if
$E\le N_s$ and its point map satisfies
$E\forces\varphi(x)\in[\dot T]\cap N_t$. Every basic pair condition admits such
an isomorphism. In $T_x$ extend $t$ by choosing the least available successor at every step.
This gives an internal branch name, so the hypothesis on $\dot T$ and
the local inverse argument above apply. Moreover every such $D$ satisfies
$D\le\pi_y(r)$. Indeed, for any nonzero basic $b\le D$, the internal branch
$\varphi(x)$ on $\varphi^{-1}(b)$ witnesses that $[T_x]\cap N_t\cap b$ is
nonempty. The second forcing can put its generic branch in that relatively
open set. Hence $r\wedge i_y(b)>0$.

\begin{lemma}\label{lem:levy-replacement}
Suppose $\varphi:E\to D$ is a local isomorphism through $r=(s,t)$,
$0<a\le E$, and $\beta<\kappa$. There are basic cylinders $U\le N_s$ and
$W\le\varphi(a)$, each determined by a string containing an entry greater
than $\beta$, and a local isomorphism through $r$ with domain contained in $U$
and range $W$, such that $\varphi^{-1}(W)\wedge U=0$.
\end{lemma}
\begin{proof}
Let $H$ be a ground target name for the inverse point map of $\varphi$ on
$D$, made total off $D$. Since $0<\varphi(a)\le\pi_y(r)$, the condition
$r\wedge i_y(\varphi(a))$ is nonzero. By Lemma~\ref{lem:levy-pair} it forces
$x\ne H[y]$. Strengthen it to decide incompatible prefixes of these two points.
Extend the input prefix past an entry greater than $\beta$, and then strengthen
to a basic pair condition $r'=(u,t')$ with this longer input prefix included
in $u$. If $v$ is the decided prefix of $H[y]$, then $u\perp v$ and
\[
 \pi_y(r')\le
 \varphi(a)\wedge \big\|v\subseteq H[\dot y]\big\|_{\BB}.
\]
Choose a local isomorphism through $r'$. Its range is below $\pi_y(r')$.
Restrict this isomorphism to a nonzero basic subcylinder $W$ of its range,
lengthening the defining string of $W$ to include an entry greater than $\beta$.
Its domain lies in $U=N_u$, whereas the old inverse image
$\varphi^{-1}(W)$ lies in $a\wedge N_v$.
\end{proof}

\subsection{Selecting a branch}\label{subsec:selecting-a-branch}

For each basic input condition $N_u\le N_{p_*}$ we shall choose disjoint basic
subconditions $U_0,U_1$ and a nonempty basic target cylinder $V$ so that
\begin{equation}\label{eq:levy-overlap}
 \begin{gathered}
 \text{for each }i<2\text{ and each nonempty basic }O\le V,\text{ some}\\
 \text{basic subcondition of }U_i\text{ forces } f(x)\in O .
 \end{gathered}
\end{equation}
To do this we build in $L$ a sequence of length $\kappa$ of permanent decisions.
A decision $(S_\delta,t_\delta)$ at stage $\delta$ specifies a basic source
cylinder $S_\delta$ forcing $t_\delta\in\dot T$ and requires
$t_\delta\subset f(x)$ there. If the new source cylinder lies inside an old
one, the new prefix must extend the old prefix.

When $(U,V)$ is registered as one of the two requirements for $u$, we attach
to it local isomorphisms $\varphi_j:E_j\to D_j$ with $E_j\le U$ and $D_j\le V$.
If $\varphi_j$ is added at stage $b_j$, its remaining domain just before
stage $\alpha>b_j$ is
\[
 A_{j,\alpha}=E_j\wedge
  \neg\bigvee\{S_\delta:b_j<\delta<\alpha\text{ and a decision was made at }\delta\}.
\]
At installation its values must be consistent with all decisions made so far.
At each $\alpha<\kappa$ we maintain
\begin{equation}\label{eq:levy-cover}
 V\le \bigvee_{j}\varphi_j(A_{j,\alpha})
\end{equation}
for every previously registered $(U,V)$, with the join over the isomorphisms
attached to that requirement.

We schedule three kinds of tasks. For each $u$ we register two requirements
$(U_0,V)$ and $(U_1,V)$ as in~\eqref{eq:levy-overlap}. For each $(u,n)$ we make a
decision of length at least $n$ below $N_u$. For each registered $(U,V)$ and each
nonempty basic $O\le V$, we make a decision below $U$ forcing $f(x)\in O$.
In $L$ there are $\kappa$ possible tasks, since their indices are finite tuples of
basic cylinders and natural numbers. Fix a schedule of length $\kappa$ in which each
possible task occurs unboundedly often. Carry out each task at its first scheduled
stage for which its prerequisites are met.

At each stage we record the basic source cylinder of any new permanent decision, a basic
source cylinder enclosing each newly installed domain, and the basic target cylinders used in
the replacements below. Only finitely many cylinders are recorded per stage. Choose
$\beta<\kappa$ bounding the entries of all strings previously recorded.
The defining string of each newly recorded cylinder must contain an entry greater than $\beta$.
Such a bound exists because each proper stage $\alpha<\kappa$ is countable
in $L$.

For use at a source step, start below any prescribed positive source
condition and choose a basic subcylinder whose string has an entry above
$\beta$. Every previously recorded source cylinder meeting it is determined by
an ancestor of that string. There are finitely many ancestors and each recorded
string occurs at only one stage. Refining finitely many more times, we
obtain a basic $N_s$ on which the applicable earlier decisions are constant
and such that, for every $j$, either $N_s\le A_{j,\alpha}$ or
$N_s\wedge A_{j,\alpha}=0$. Only finitely many $j$ give the first alternative.
Indeed, only finitely many enclosing cylinders and decision cylinders enter the
test, and we can decide each relevant Boolean domain on a stronger basic
condition. Let $t$ be the longest applicable output prefix, or the empty
string if there is none. The pair $r=(s,t)$ is a condition. Every isomorphism
still defined throughout $N_s$ restricts there to an isomorphism through $r$.

To register the two requirements for $u$, choose $N_s\le N_u$ as above and take a
local isomorphism through $r$. Apply Lemma~\ref{lem:levy-replacement} with $a$ its domain.
The lemma gives a second isomorphism whose source lies in a cylinder $U_1$
disjoint from the first inverse image of $W$. Choose a basic $U_0$ inside that
inverse image with an entry above $\beta$ in its string. Restrict the first
map to $U_0$, and then restrict both maps to a common nonempty basic target
$V$ inside the image of $U_0$. Install the resulting maps for $(U_0,V)$ and
$(U_1,V)$ respectively. Their ranges are all of $V$. This step makes no
permanent decision.

To make a permanent decision, we first choose a source cylinder $N_s$ and a prefix
$t'$ according to the task. If the task specifies a basic $O\le V$ for a registered
requirement $(U,V)$, \eqref{eq:levy-cover} supplies an isomorphism whose
remaining target image meets $O$. Start the preceding source refinement
below the inverse image of this intersection. Then $N_s\le U$, and its old
prefix $t$ is compatible with the string defining $O$, since the isomorphism witnesses both
on $N_s$. Let $t'$ be the longer string. If the task specifies $(u,n)$, carry out
the source refinement below $N_u$. Since $N_s$ forces $t$ to belong to the pruned tree,
we can choose an extension $t'$ of $t$ of length at least $n$ and strengthen $N_s$
to force $t'\in\dot T$. In either case, the resulting $N_s$ forces $t'\in\dot T$.

Write $\varphi_1,\ldots,\varphi_m$ for the isomorphisms whose remaining
domains contain $N_s$, restricting them temporarily to $N_s$. Start with
$a=N_s$. For $j=1,\ldots,m$, apply Lemma~\ref{lem:levy-replacement} to
$\varphi_j$ through $(s,t)$ and the current $a$. Keep the replacement
isomorphism, its enclosing source cylinder $U'_j$ and its target $W_j$,
and replace $a$ by $\varphi_j^{-1}(W_j)$. These are successive positive
refinements. Choose a basic cylinder $S\le a$ with an entry above $\beta$
in its string, and put $(S,t')$ into the permanent decisions.

For every $j$ the entire lost image $\varphi_j(S)$ lies in $W_j$, and
$S\wedge U'_j=0$. Install the replacement with the same requirement as
$\varphi_j$. Its domain is inside $N_s$, hence inside the required source
cylinder, and its outputs extend the old prefix $t$. Its domain misses
the new decision. Thus it is consistent with all decisions and covers the
entire loss. Its source can overlap the sources of other replacements.
Record $W_j$ as the target cylinder for this replacement. All other
previous maps were already absent from $S$. This proves
\eqref{eq:levy-cover} at the successor stage.

It remains to verify~\eqref{eq:levy-cover} at a limit stage $\lambda<\kappa$.
Fix a previously registered $(U,V)$ and a bound $\beta<\kappa$ on every entry recorded
before $\lambda$. The basic $N_d\le V$ whose defining strings have an entry above
$\beta$ have join $V$. Fix such a $d$.
Every earlier replacement target meeting $N_d$ is determined by an ancestor
of $d$. Choose $\gamma<\lambda$ after registration and after all of those
finitely many replacement stages. For each isomorphism $j$ of this requirement
already installed before $\gamma$, set
\[
 B_{j,\xi}=N_d\wedge\varphi_j(A_{j,\xi})
 \qquad (\gamma\le\xi\le\lambda).
\]
The sequence is decreasing. No successor step after stage $\gamma$ can remove
a positive part from $B_{j,\xi}$, since its recorded replacement target would
then meet $N_d$. At a limit $\eta$ the definition of $A_{j,\eta}$ and completeness
of $\varphi_j$ give
\[
 B_{j,\eta}=\bigwedge_{\gamma\le\xi<\eta}B_{j,\xi}.
\]
It follows that $B_{j,\lambda}=B_{j,\gamma}$. Coverage at $\gamma$ gives
$N_d\le\bigvee_j\varphi_j(A_{j,\lambda})$. Taking the join over $d$ proves
\eqref{eq:levy-cover} at $\lambda$.

This completes the recursion in $L$. For a generic input $x$, the union of the permanent
prefixes applying to it is a branch $f(x)$ through $T_x$. Indeed, for each $n$ the decisions of
length at least $n$ are dense. Whenever a new decision cylinder meets an old
one, its entry above the previous bound prevents it from being an ancestor
of the old cylinder. Thus the applicable prefixes form a chain. The output
is generic by the hypothesis on $T_x$, and \eqref{eq:levy-overlap} gives
$V\le\pi_{i_f}(U_0)\wedge\pi_{i_f}(U_1)$ for each registered pair.

Since $p_*$ forces every branch of $T_x$ to generate the full extension, on some basic source condition
$N_a$ the map induces an isomorphism onto a target condition. The target
support of any $e\le N_a$ is unchanged by restricting the homomorphism to
$N_a$, since $e\le i_f(b)$ is equivalent to $e\le N_a\wedge i_f(b)$.
The disjoint $U_0,U_1$ registered below $N_a$ must therefore have disjoint
target supports, a contradiction. This proves Theorem~\ref{thm:single-collapse}.

\begin{corollary}\label{cor:single-od}
In the single $\PP$-generic extension of $L$, $\mathrm{PK}(\mathrm{OD})$ fails.
\end{corollary}
\begin{proof}
Using the canonical real coding $J$ fixed above, in $L[g]$ the set
\[
 A=\{J(h):h\in\kappa^\omega\text{ is $\PP$-generic over }L
                        \text{ and }V=L[h]\}
\]
is ordinal-definable. Partition $\kappa$ in $L$ into $\kappa$ pairs.
Independently swapping or fixing each pair gives $(2^\kappa)^L$
ground permutations $\rho$ of $\kappa$. The points $\rho\circ g$ are
distinct full generators, since $g$ is onto and each $\rho$ acts by an automorphism of $\PP$.
Hence $A$ is uncountable, since forcing of size $\kappa$ preserves
$(\kappa^+)^L$.

In $L[g]$ the space $\kappa^\omega$ is Polish and $J$ is continuous. If $K$
were a copy of Cantor space in $A$, then $J^{-1}[K]$ would be an uncountable
closed subset of $\kappa^\omega$. By the perfect set theorem for closed sets
it would contain a copy of Cantor space, contrary to
Theorem~\ref{thm:single-collapse}.
\end{proof}

Since the generic filter is determined by its union $g$, the single L\'evy extension is
$L[J(g)]$. This completes the proof of Theorem~\ref{thm:summary}.

\begin{remark}
The construction of Theorem~\ref{thm:single-collapse} also applies to Cohen forcing,
presented as $\omega^{<\omega}$. The original recursion had $\omega_1^L$ tasks in $L$.
With $\kappa=\omega$, the strings indexing registrations and decisions belong to the
countable set $\omega^{<\omega}$, so in $L$ we can schedule every task infinitely
often in $\omega$ stages. At each stage only finitely many strings have been recorded,
so their entries have a finite bound. Every Cohen generic has arbitrarily large entries.
Lemmas~\ref{lem:levy-pair} and~\ref{lem:levy-replacement} and the successor construction
therefore apply, and there are no nonzero limit stages below $\omega$. Thus, if $c$ is
Cohen-generic over $L$, the OD set of Cohen reals $h$ with $L[h]=L[c]$ contains no Cantor
set in $L[c]$. Since it contains $c$ and is disjoint
from $L$, this provides a negative answer to Schilhan's Question~44 about $\mathrm{MS}(\mathrm{OD})$
\cite{Schilhan}.
\end{remark}

\section{AI Disclosure}
The proof Theorem~\ref{thm:single-collapse} was obtained by the author while studying the (entirely AI-generated) proofs in \cite{Chen26}.

\end{document}